\documentclass[leqno, fleqn]{amsart}%
\usepackage{amsmath}
\usepackage{amsfonts}
\usepackage{amssymb}%

\makeatletter \DeclareMathSymbol{\Gamma}{\mathalpha}{letters}{"00}
\DeclareMathSymbol{\Delta}{\mathalpha}{letters}{"01}
\DeclareMathSymbol{\Theta}{\mathalpha}{letters}{"02}
\DeclareMathSymbol{\Lambda}{\mathalpha}{letters}{"03}
\DeclareMathSymbol{\Xi}{\mathalpha}{letters}{"04}
\DeclareMathSymbol{\Pi}{\mathalpha}{letters}{"05}
\DeclareMathSymbol{\Sigma}{\mathalpha}{letters}{"06}
\DeclareMathSymbol{\Upsilon}{\mathalpha}{letters}{"07}
\DeclareMathSymbol{\Phi}{\mathalpha}{letters}{"08}
\DeclareMathSymbol{\Psi}{\mathalpha}{letters}{"09}
\DeclareMathSymbol{\Omega}{\mathalpha}{letters}{"0A}
\DeclareMathSymbol{\varGamma}{\mathalpha}{operators}{"00}
\DeclareMathSymbol{\varDelta}{\mathalpha}{operators}{"01}
\DeclareMathSymbol{\varTheta}{\mathalpha}{operators}{"02}
\DeclareMathSymbol{\varLambda}{\mathalpha}{operators}{"03}
\DeclareMathSymbol{\varXi}{\mathalpha}{operators}{"04}
\DeclareMathSymbol{\varPi}{\mathalpha}{operators}{"05}
\DeclareMathSymbol{\varSigma}{\mathalpha}{operators}{"06}
\DeclareMathSymbol{\varUpsilon}{\mathalpha}{operators}{"07}
\DeclareMathSymbol{\varPhi}{\mathalpha}{operators}{"08}
\DeclareMathSymbol{\varPsi}{\mathalpha}{operators}{"09}
\DeclareMathSymbol{\varOmega}{\mathalpha}{operators}{"0A}

\newcommand{\allmodesymb}[2]{\relax\ifmmode{\mathchoice
{\mbox{\fontsize{\tf@size}{\tf@size}#1{#2}}}
{\mbox{\fontsize{\tf@size}{\tf@size}#1{#2}}}
{\mbox{\fontsize{\sf@size}{\sf@size}#1{#2}}}
{\mbox{\fontsize{\ssf@size}{\ssf@size}#1{#2}}}} \else
\mbox{#1{#2}}\fi}

\makeatother
\makeatletter
\renewcommand*\subjclass[2][2000]{%
  \def\@subjclass{#2}%
  \@ifundefined{subjclassname@#1}{%
    \ClassWarning{\@classname}{Unknown edition (#1) of Mathematics%
      Subject Classification; using '2000'.}%
  }{%
    \@xp\let\@xp\subjclassname\csname subjclassname@#1\endcsname%
  }%
} \makeatother

\theoremstyle{plain}

\theoremstyle{remark}

\allowdisplaybreaks \numberwithin{equation}{section}
\begin{document}
\title{The equivalent classical metrics on the Cartan-Hartogs Domains}
\date{18th September 2005}

\author{Weiping YIN }
\address{W. YIN: Dept. of Math., Capital Normal Univ., Beijing 100037, China}
\email{wyin@mail.cnu.edu.cn; wpyin@263.net}
\author{An WANG }
\address{An WANG: Dept. of Math., Capital Normal Univ., Beijing
100037, China} \email{wangan@mail.cnu.edu.cn} \subjclass{32H15,
32F07, 32F15} 
\thanks{Project supported in part by NSF of China (Grant NO.
10471097) and the. Doctoral Programme Foundation of NEM of China.}

\begin{abstract}
In this paper we study the complete invariant metrics on
Cartan-Hartogs domains which are the special types of Hua domains.
Firstly, we introduce a class of new complete invariant metrics on
these domains, and prove that these metrics are equivalent to the
Bergman metric. Secondly, the Ricci curvatures under these new
metrics are bounded from above and below by the negative
constants. Thirdly, we estimate the holomorphic sectional
curvatures of the new metrics, we prove that the holomorphic
sectional curvatures are bounded from above and below by the
negative constants. Finally, by using these new metrics and Yau's
Schwarz lemma we prove that the Bergman metric is equivalent to
the Einstein-K\"ahler metric. That means the Yau's conjecture is
true on Cartan-Hartogs domain.

\end{abstract}
\baselineskip0.6cm \maketitle

The concept of Hua domain was introduced by Weiping Yin in 1998.
Since then, many good results have been obtained. The Bergman kernel
functions are given in explicit forms[1-17]. The comparison theorems
for Bergman metric and Kobayashi metric are proved on Cartan-Hartogs
domains[18-21]. The explicit form of the Einstein-K\"{a}hler metric
is got on non-symmetric domain which is the first time in the
world[22-26], etc. In this paper we will study the equivalence
between the classical metrics. There are many deep results on this
subject. Let $\omega_B(D)$, $\omega_C(D)$, $\omega_K(D)$,
$\omega_{EK}(D)$ be the Bergman metric, Carath\'{e}odory metric,
Kobayashi metric and Einstein-K\"{a}hler metric on bounded domain
$D$ in $\mathbb{C}^n$ respectively. Then we have
$\omega_C(D)\leqslant 2\omega_B(D)$[27,28], $\omega_C(D)\leqslant
\omega_K(D)$[29], $\omega_C(D)=\omega_K(D)$ if $D$ is the convex
domain[30], $\omega_B(D)=\omega_{EK}(D)$ if $D$ is the bounded
homogeneous domain in $\mathbb{C}^n$[31,P.300]. For the
$\omega_B(D)$ and $\omega_K(D)$, no relationship is known. People
had hoped that the inequality $\omega_B(D)\leqslant C \omega_K(D)$
for some universal constant $C$ would hold, but in 1980 Diederich
and Foraess [32] showed that there exist pseudoconvex domain in
$\mathbb{C}^3$ where the quotient $\omega_B(D)/\omega_K(D)$ is
unbounded. If the inequality $\omega_B(D)\leqslant C \omega_K(D)$
holds then we say that the comparison theorem for Bergman metric and
Kobayashi metric on $D$ holds. The above comparison theorem [18-21]
are in this sense. Recently, Kefeng Liu, Xiaofeng Sun and Shing-Tung
Yau study the equivalence between the classical metrics on
Teichm\"{u}ller spaces and moduli spaces[33-35]. They proved that on
Teichm\"{u}ller spaces and moduli spaces the four classical metrics
$\omega_B(D), \omega_C(D), \omega_K(D), \omega_{EK}(D)$ are
equivalent. Especially, they proved the conjectures of Yau about the
equivalence between the Einstein-K\"{a}hler metric and the
Teichm\"{u}ller metric and also its equivalence with the Bergman
metric.

In this paper we study the complete invariant metrics on
Cartan-Hartogs which is the special types of Hua domains. We prove
that the Yau's conjecture is true on the Cartan-Hartogs domains.
This paper is organized as follows. We introduce a class of new
complete invariant metrics on Cartan-Hartogs and prove that these
metrics are equivalent to the Bergman metrics of the
Cartan-Hartogs domains in the first section. In the second section
we prove that the Ricci curvatures of these new metrics are
bounded from above and below by the negative constants. We prove
also that the holomorphic sectional curvatures of these new
metrics on Cartan-Hartogs domains are bounded from above and below
by the negative constants in the third section. In the fourth
section, by using these new metrics and Yau's Schwarz lemma[36] we
prove that the Bergman metric is equivalent to the
Einstein-K\"{a}hler metric. That means the Yau's conjecture[37] is
true on Cartan-Hartogs domain.

The Cartan-Hartogs domains are defined as follows:
$$Y_I=\{ W\in \mathbb{C}^N,Z\in R_I(m,n):|W|^{2K} <\det(I-Z\overline{Z}^t),K>0\},$$
$$Y_{II}=\{ W\in \mathbb{C}^N,Z\in R_{II}(p):|W|^{2K}<\det(I-Z\overline{Z}^t),K>0 \},$$
$$Y_{III}=\{ W\in \mathbb{C}^N,Z\in R_{III}(q):|W|^{2K}<\det(I-Z\overline{Z}^t),K>0 \},$$
$$Y_{IV}=\{ W\in \mathbb{C}^N,Z\in R_{IV}(n):|W|^{2K}<(1-2Z\overline{Z}^t+|ZZ^t|^2), K>0 \}.
$$
Where $\det$ is the abbreviation of determinant; $Z^t$ indicates
the transpose of $Z$, $\overline{Z}$ denotes the conjugate of $Z$;
and the $R_I(m,n)$, $R_{II}(p)$, $R_{III}(q)$, $R_{IV}(n)$ denote
the classical domains in the sense of Hua[38](they are also called
Cartan domains).

\section {New complete invariant metrics on $Y_I$}
In this section we will introduce the new complete invariant
metrics on $Y_I$, and prove that these new metrics are equivalent
to the Bergman metric on $Y_I$.

1.1. New invariant metrics

1.1.1. Suppose $(Z,W)\in Y_I, Z=(z_{ij}), W=(w_1,w_2,\dots, w_N)$
and let
$$Z_1=(z_1,z_2,\dots, z_{mn})=(z_{11},z_{12},\dots, z_{1n},
\dots, z_{m1},z_{m2},\dots, z_{mn}),$$
$$Z_2=(z_{mn+1},z_{mn+2},\dots, z_{mn+N})=(w_1,w_2,\dots, w_N),$$
then the point $(Z,W)$ can be denoted by a vector $z$ with $mn+N$
entries, that is
$$z=(Z_1,Z_2)=(z_1,z_2,\dots, z_{mn}, z_{mn+1},z_{mn+2},\dots, z_{mn+N}).$$

1.1.2. Let
$$
G_{\lambda}=G_{\lambda}(Z,W)=Y^{\lambda}[\det(I-Z\overline
Z^t)]^{-(m+n+\frac{N}{K})}, \lambda > 0,
$$
$$T_{\lambda I}(Z,W;\overline{Z},\overline{W})
=(g_{i\overline{j}})=\left(\frac{\partial^2\log G_{\lambda}}
{\partial z_i\partial \overline{z}_{j}}\right),$$
where
$$Y=(1-X)^{-1}, X=|W|^2[\det(I-Z\overline Z^t)]^{-\frac 1K}.$$

1.1.3. The following mappings are the holomorphic automorphism of
$Y_I$, which map the point $(Z,W)\in Y_I$ onto point $(0,W^*)$:
$$ \left\{
\begin{array}{lll}
W^*& =& W\det{(I-Z_0\overline Z_0^t)^{\frac1{2K}}}
\det{(I-Z\overline Z_0^t)^{-\frac1K}} \\
Z^*& =& A(Z-Z_0){(I-{\overline
Z_0^t}Z)^{-1}}{{\overline{D}}\!^{-1}}
\end{array}\right.$$
Where $\overline{A}^tA=(I-Z_0\overline{Z}_0^t)^{-1},
\overline{D}^tD=(I-\overline{Z}_0^tZ_0)^{-1}.$ The set of these
mappings is denoted by ${\rm Aut}(Y_I).$

1.1.4. Let $f\in{{\rm Aut}(Y_I)}$ and $Z_0=Z$, then one has

$$ T_{\lambda I}(Z,W;\overline{Z},\overline{W})=J_f|_{Z_0=Z}
T_{\lambda I}(0,W^*;0,\overline{W^*}) \overline{J_f}^t|_{Z_0=Z},$$
the Jacobian matrix of $f$ is equal to
$$J_f=\left(
\begin{array}{cc}
\frac{\partial Z^*}{\partial Z_1}&\frac{\partial W^*}{\partial
Z_1}
\vspace{2mm}\\
\frac{\partial Z^*}{\partial Z_2}&\frac{\partial W^*}{\partial
Z_2}
\end{array}\right)$$

Let $$ J_f|_{Z_0=Z}=\left(
\begin{array}{cc}
J_{11}&J_{12}\\
J_{21}&J_{22}
\end{array}
\right).
$$
Then one has $$
\begin{array}{ll}
J_{11}=A^t\cdot\!\!\times \overline{D}^t,\\
\\
J_{12}=\frac{1}{K}\det(I-Z\overline{Z}^t)^{-\frac{1}{2K}}E(Z)^tW,\\
\\
J_{21}=0,\\
\\
J_{22}=\det(I-Z\overline{Z}^t)^{-\frac{1}{2K}}I.
\end{array}$$
where
$$E(Z)=\left(
{\rm tr}[(I-Z\overline{Z}^t)^{-1} I_{11}\overline{Z}^t], {\rm
tr}[(I-Z\overline{Z}^t)^{-1} I_{12}\overline{Z}^t], \cdots, {\rm
tr}[(I-Z\overline{Z}^t)^{-1} I_{mn}\overline{Z}^t] \right)$$ is
the $1\times mn$ matrix. And $I_{\alpha \beta}$ is $m\times n$
matrix, the $(\alpha \beta)$-th entry of  $I_{\alpha \beta}$ ,
i.e. the entry located at the junction of the $\alpha$-th row and
$\beta$-th column of $I_{\alpha \beta}$ is 1, and its others
entries are zero. The meaning of $ [A \cdot\!\!\times A] $ can be
found in [39] or in 1.1.6 below.

1.1.5. By computations,one has
$$T_{\lambda
I}(0,W^*;0,\overline{W^*})=\left(\begin{array}{cc}
(\frac{\lambda}{K} YX+m+n+\frac{N}{K})I&0\\
0& \lambda YI+\lambda Y^2{\overline {W^*}^t}W^*
\end{array}\right).$$
Where $Y, X$ see 1.1.2.

1.1.6. Let
$$
T_{\lambda I}(Z,W;\overline{Z},\overline{W})
=J_f|_{Z_0=Z}T_{\lambda
I}(0,W^*;0,\overline{W^*})\overline{J_f}^t|_{Z_0=Z}
=\left(\begin{array}{cc}
T_{11}&T_{12}\\
T_{21}&T_{22}
\end{array}
\right),
$$
one has
$$
\begin{array}{ll}
T_{11}=(\frac\lambda{K} YX+m+n+\frac{N}{K})(A^t\overline A
\cdot\!\!\times \overline D^tD)
+\frac{\lambda}{K^2} Y^2XE(Z)^t\overline{E(Z)},\\
\\
T_{12}=\frac{1}{K}\det(I-Z\overline{Z}^t)^{-\frac{1}{K}}\lambda
Y^2E(Z)^tW,\\
\\
T_{21}=\overline T^t_{12},\\
\\
T_{22}=\lambda Y\det(I-Z\overline{Z}^t)^{-\frac{1}{K}}I+
\det(I-Z\overline{Z}^t)^{-\frac{2}{K}}\lambda Y^2\overline W^tW.
\end{array}$$
The computations form 1.1.3 to 1.1.6 can be found in [19] for
detail. The definition of $\cdot\!\!\times$ is the following (see
[39]).

Definition: For the $r\times s$ matrix
$$A=\left(\begin{array}{ccc}
a_{11}&\dots & a_{1s}\\
\dots& \dots & \dots \\
a_{r1} &\dots & a_{rs}
\end{array}\right).$$
and $p\times q$ matrix
$$B=\left(\begin{array}{ccc}
b_{11}&\dots & b_{1q}\\
\dots& \dots & \dots \\
b_{p1} &\dots & b_{pq}
\end{array}\right).$$ The $\cdot\!\!\times $ of these two matrices
is as the follows:
$$A \cdot\!\!\times B=\left(\begin{array}{ccc}
a_{11}B&\dots & a_{1s}B\\
\dots& \dots & \dots \\
a_{r1}B &\dots & a_{rs}B
\end{array}\right).$$
Which is a $rp\times sq$ matrix.

1.1.7. From 1.1.5 and 1.1.6, one has $T_{\lambda
I}(Z,W;\overline{Z},\overline{W})>0$, and then from the definition
of $G_{\lambda}$ in 1.1.2, the $G_{\lambda}$ generates an invariant
metric $\omega_{G_{\lambda}}(Y_I)$ of $Y_I$.

1.2. New invariant metrics are equivalent to Bergman metric

Definition: Two metrics  $\mathcal{B}$ and $\mathcal{E}$ of domain
$\Omega$ in $\mathbb{C}^n$ are called equivalent, if they are
quasi-isometric to each other in the sense that
$$ b\leqslant \frac{\mathcal{B}}{\mathcal{E}}\leqslant a.$$
for two positive
constants $a$ and $b$. We will write this as $\mathcal{B}\sim
\mathcal{E}.$

The Bergman kernel function of $Y_I$ has the following form[03]:
$$K_{Y_I}=K^{-mn}{\pi}^{-mn+N}G(X)\det(I-Z\overline{Z}^t)^{- (m+n+\frac{N}{K})}.$$
Where $\displaystyle
G(X)=\sum_{j=0}^{mn+1}b_j\Gamma(N+j)(1-X)^{-(N+j)}$, let
$$
\begin{array}{ll}
P(x)=&(x+1)[(x+1+Kn)(x+1+K(n-1))\dots(x+1+K)]\\[3mm]
&[(x+1+K(n+1))(x+1+Kn)\dots(x+1+2K)]\\[3mm]
&[(x+1+K(n+2))(x+1+K(n+1))\dots(x+1+3K)]\dots\\[3mm]
&[(x+1+K(n+m-1))(x+1+K(n+m-2))\dots(x+1+mK)].
\end{array}
$$
Then $b_0=P(-1)=0.$ And the other $b_j(j=1, 2, \cdots, mn+1)$ is
determined by the following form:
$$
b_j=
\frac{P(-j-1)-\displaystyle\sum_{k=0}^{j-1}b_k(-1)^k\frac{\Gamma(j+1)}{\Gamma(j-k+1)}}{
(-1)^j\Gamma(j+1)},
$$
The Bergman kernel function of $Y_I$ generates the Bergman metric
$\omega_B(Y_I).$

1.2.1. By calculations, the metric matrix of Bergman metric
$$
T_{BI}(Z,W;\overline{Z},\overline{W}) =J_f|_{Z_0=Z}T_{BI}(Z^*,W^*;
\overline{Z^*},\overline{W^*})|_{Z^*=0}
 \overline{J_f}^t |_{Z_0=Z}$$
$$=J_f|_{Z_0=Z}\left(
\begin{array}{cc}
[\frac{1}{K}M'X+m+n+\frac{N}{K}]I&0\\
0&M'I+M''{\overline {W^*}}^t W^*
\end{array}\right) \overline{J_f}^t |_{Z_0=Z}.$$
Where
$$\log G(X)=M,\ \ \ \ \frac{\partial \log G(X)}{\partial X}=M',
\ \ \ \ \frac{\partial^2 \log G(X)}{\partial X^2} =M''.$$ And
$J_f|_{Z_0=Z}$ is same as that in 1.1.4. For the details please
see [03]. Hence
$$(\omega_B(Y_I))^2=dzJ_f|_{Z_0=Z}T_{BI}(0,W^*;
0,\overline{W^*}) \overline{J_f}^t |_{Z_0=Z}\overline{dz}^t$$
$$=dzJ_f|_{Z_0=Z}\left(
\begin{array}{cc}
[\frac{1}{K}M'X+m+n+\frac{N}{K}]I&0\\
0&M'I+M''{\overline {W^*}}^t W^*
\end{array}\right) \overline{J_f}^t |_{Z_0=Z}\overline{dz}^t.$$

1.2.2. Due to 1.1.5., one has
$$(\omega_{G_{\lambda}}(Y_I))^2=dzJ_f|_{Z_0=Z}
T_{\lambda I}(0,W^*;0,\overline{W^*})
\overline{J_f}^t|_{Z_0=Z}\overline{dz}^t
$$
$$
=dzJ_f|_{Z_0=Z}\left(
\begin{array}{cc}
[\frac{1}{K}\lambda YX+m+n+\frac{N}{K}]I&0\\
0& \lambda YI+\lambda Y^2{\overline {W^*}^t}W^*
\end{array}\right)\overline{J_f}^t|_{Z_0=Z}\overline{dz}^t
.$$ Where $Y, X$ are same as that in 1.1.2.

1.2.3. Let $dzJ_f|_{Z_0=Z}=(d\tilde{z}, d\Im)$, one has

$(\omega_{G_{\lambda}}(Y_I))^2 =(\frac{\lambda
YX}{K}+m+n+\frac{N}{K})|d\tilde{z}|^2+d\Im(\lambda YI+\lambda Y^2
\overline {W^*}^tW^*) \overline{d\Im}^t,$

($\omega_B(Y_I))^2=(\frac{1}{K}M'X+m+n+\frac{N}{K})|d\tilde{z}|^2
+d\Im(M'I+M''{\overline {W^*}}^t W^*) \overline{d\Im}^t.$

1.2.4. From [39], the vector $W^*=(w^*_1,w^*_2,\cdots,w^*_N)$ can
be written as
$$W^*=e^{i\theta}(\mu,0,\cdots,0)U,\quad \mu\geqslant 0,$$
where $U$ is the unitary matrix, hence
$$
\overline{W^*}^tW^*=\overline{U}^t\left(\begin{array}{cc}\mu&0\\
0&\mathbf{0} \end{array}\right)U,$$
$$\quad M'I^{(N)}+M''{\overline
W^{*}}^{t}W^{*}=\overline{U}^t\left(\begin{array}{cc}M'+M''\mu^2&0\\
0&M'I^{(N-1)} \end{array}\right)U,$$
$$
\lambda YI^{(N)} + \lambda Y^{2}\overline {W^*}^tW^*=
\overline{U}^t\left(\begin{array}{cc}\lambda Y+\lambda Y^2\mu^2&0\\
0&\lambda YI^{(N-1)} \end{array}\right)U>0.
$$

1.2.5. Let $d\Im\overline{U}^t=(dv,d\tilde{W})$, then
$\omega_{G_{\lambda}}(Y_I)$ and $\omega_B(Y_I)$ can be written as:
$$\begin{array}{l}
(\omega_{G_{\lambda}}(Y_I))^2=(\frac{\lambda
YX}{K}+m+n+\frac{N}{K})|d\tilde{z}|^2+
(\lambda Y+\lambda Y^2\mu^2)|dv|^2+\lambda Y|d\tilde{W}|^2,\\
(\omega_B(Y_I))^2=(\frac{1}{K}M'X+m+n+\frac{N}{K})|d\tilde{z}|^2+
(M'+M''\mu^2)|dv|^2+M'|d\tilde{W}|^2.
\end{array}$$

1.2.6. $T_{\lambda I}(Z,W;\overline{Z},\overline{W})>0$ and
$T_{BI}(Z,W;\overline Z,\overline{W})>0$, these imply
$$\frac{1}{K}M'X+m+n+\frac{N}{K}>0,\quad M'+M''\mu^2 >0, \quad  M' >0,\quad Y >0.$$

1.2.7. Let
$$
\Phi(X)=\frac{\frac{1}{K}M'X+m+n+\frac{N}{K}}{\frac{\lambda
YX}{K}+m+n+\frac{N}{K}},\quad \Psi(X)=\frac{M'+M''\mu^2}{\lambda
Y+\lambda Y^2\mu^2},\quad \Upsilon(X)=\frac{M'}{\lambda Y},
$$
then all of $\Phi(X),\Psi(X)$, $\Upsilon(X)$ is positive continues
function of $X$ on the interval $[0,1)$. If
$$ \lim_{X\rightarrow 1}\Phi(X),\ \ \ \ \lim_{X\rightarrow 1}\Psi(X),
\ \ \ \ \lim_{X\rightarrow 1}\Upsilon(X)$$ are existent and
positive, then all of $\Phi(X), \Psi(X), \Upsilon(X)$ have the
positive maximum and the positive minimum on $[0, 1)$.

1.2.8. Because
$$G(X)=\sum_{j=0}^{mn+1}b_j\Gamma(N+j)(1-X)^{-(N+j)}=
\sum_{j=0}^{mn+1}b_j\Gamma(N+j)Y^{(N+j)},
$$
therefore
$$
\frac{dG(Y)}{dX}=G'(X)=\sum_{j=0}^{mn+1}b_j\Gamma(N+j+1)Y^{(N+j+1)},$$
$$\quad
\frac{d^2G(Y)}{dX^2}=G''(X)=\sum_{j=0}^{mn+1}b_j\Gamma(N+j+2)Y^{(N+j+2)}.
$$ Let
$$M'=G'(X)G^{-1}(X),$$ $$M''=G''(X)G^{-1}(X)-G'(X)^2G^{-2}(X).$$

1.2.9. We will compute the limits of $\Phi(X)$, $\Psi(X)$ and
$\Upsilon(X)$ when $X$ tends to 1.

By calculations, one has
$$
\lim_{X\rightarrow
1}\Phi(X)=\lim_{Y\rightarrow\infty}\Phi(X)=\frac{mn+N+1}{\lambda}.
$$
Hence there exists $0<\nu<\delta$ such that $$0<\nu\leqslant
\Phi(X)\leqslant \delta.
$$

1.2.10. Similarly, one has
$$
\lim_{X\rightarrow
1}\Psi(X)=\lim_{Y\rightarrow\infty}\Psi(X)=\frac{mn+N+1}{\lambda},
$$
and
$$
\lim_{X\rightarrow1}\Upsilon(X)=\lim_{Y\rightarrow\infty}\Upsilon(X)
=\frac{mn+N+1}{\lambda}.
$$
Therefore there exist $\zeta,\ \eta,\ \rho$ and $\varrho$ such
that
$$
0<\zeta\leqslant\Psi(X)\leqslant\eta,\quad
0<\rho\leqslant\Upsilon(X)\leqslant\varrho.
$$

1.2.11. Let $a^2=\max\{\mu,\eta,\varrho\}$ and
$b^2=\min\{\nu,\zeta,\rho\}$, then one has $0< b\leqslant
\frac{\omega_B(Y_I)}{\omega_{G_{\lambda}}(Y_I)}\leqslant a.$
Therefore the following theorem is proved

Theorem: All of the above new complete invariant metrics is
equivalent to the Bergman metric on $Y_{I}$. That is
$$\omega_B(Y_I)\sim \omega_{G_{\lambda}}(Y_I).$$
Because the Bergman metric $\omega_B(Y_I)$ is complete ([40]), hence
the new metric $\omega_{G_{\lambda}}(Y_I)$ is also complete.

1.2.12. By using the same idea and method. If we introduce the
following functions
$$G_{\lambda}=Y^{\lambda}[\det(I-Z\overline Z^t)]^{-(p+1+\frac{N}{K})}, \lambda>0;$$
$$(Y=(1-X)^{-1}, X=|W|^2[\det(I-Z\overline Z^t)]^{-\frac 1K}, (Z,W)\in{Y_{II}})$$
$$G_{\lambda}=Y^{\lambda}[\det(I-Z\overline Z^t)]^{-(q-1+\frac{N}{K})}, \lambda > 0;$$
$$(Y=(1-X)^{-1}, X=|W|^2[\det(I-Z\overline Z^t)]^{-\frac 1K}, (Z,W)\in{Y_{III}})$$
$$G_{\lambda}=Y^{\lambda}\beta(Z,\overline{Z})^{-(n+\frac{N}{K})}, \lambda > 0;$$
$$(Y=(1-X)^{-1}, X=|W|^2[\beta(Z,\overline{Z})]^{-\frac 1K},
\beta(Z,\overline{Z})=1+|ZZ^t|^2-2Z\overline{Z}^t,
(Z,W)\in{Y_{IV}})$$ for $Y_{II}, Y_{III}, Y_{IV}$ respectively,
then above functions generate the complete invariant metrics and
equivalent to the Bergman metrics on $Y_{II}, Y_{III}, Y_{IV}$
respectively.

\section {Ricci curvatures of new complete invariant metrics}
In this section we will prove that the Ricci curvatures of new
complete invariant metrics are bounded from above and below by the
negative constants on the Cartan-Hartogs domains.

2.1. Due to the definition, the Ricci curvature $Ric_{\lambda I}$ of
$\omega_{G_{\lambda}}(Y_I)$ on $Y_I$ has the following form
$$Ric_{\lambda I}=-\frac{dz(\frac{\partial^2
\log[\det T_{\lambda I}(Z,W;\overline{Z},\overline{W})]}{\partial
z_i\partial \overline{z}_j})\overline{dz}^t}{dzT_{\lambda
I}(Z,W;\overline{Z},\overline{W})\overline{dz}^t}.$$

2.2. Let $$G_{I}(X)= \left(\frac{\lambda
Y}{K}+m+n+\frac{N-\lambda}{K}\right)^{mn}\lambda^NY^{N+1},$$ then
from 1.1.4 and 1.1.5., one has
$$\det T_{\lambda I}(Z,W;\overline{Z},\overline{W})=G_{I}(X)[\det(I-Z\overline{Z}^t)]^
{-(m+n+\frac{N}{K})}.$$

2.3. Similar to 1.2.1, one has
$$dz\left(\frac{\partial^2\log[\det T_{\lambda I}
(Z,W;\overline{Z},\overline{W})]}{\partial z_i\partial
\overline{z}_j}\right)\overline{dz}^t:=(\omega_{det}(Y_I))^2$$
$$=dzJ_f|_{Z_0=Z}\left(
\begin{array}{cc}
[\frac{1}{K}M'_IX+m+n+\frac{N}{K}]I&0\\
0&M'_I I+M_I''{\overline {W^*}}^t W^*
\end{array}\right) \overline{J_f}^t |_{Z_0=Z}\overline{dz}^t.$$
Where $$\log G_I(X)=M_I,\ \ \ \frac{\partial \log G_I(X)}{\partial
X}=M'_I,\ \ \ \frac{\partial^2 \log G_I(X)}{\partial X^2}=M''_I.
$$ Here $J_f|_{Z_0=Z}$ is same as that in 1.1.4. In the following we will
prove $$\left(\frac{\partial^2\log[\det T_{\lambda
I}(Z,W;\overline{Z},\overline{W})]}{\partial z_i\partial
\overline{z}_j}\right)>0.$$ That is the
$$\left(
\begin{array}{cc}
[\frac{1}{K}M'_IX+m+n+\frac{N}{K}]I&0\\
0&M'_I I+M_I''{\overline {W^*}}^t W^*
\end{array}\right)$$ is positive definite matrix.
By calculations one has
$$M'_I=\frac{mn\lambda Y^2}{\lambda Y+K(m+n)+N-\lambda}+(N+1)Y>0,$$
hence
$$\frac{1}{K}M'_IX+m+n+\frac{N}{K}=\frac{mn\lambda Y(Y-1)}{K(\lambda Y+K(m+n)+N-\lambda)}$$
$$+\frac{(N+1)(Y-1)}{K}+m+n+\frac{N}{K}>0.$$
By calculations, one has
$$M''_I=\frac{mn\lambda Y^3(\lambda Y+2Km+2Kn+2N-2\lambda)}{(\lambda Y+K(m+n)+N-\lambda)^2}
+(N+1)Y^2.$$ Because $W^*=(w^*_1,w^*_2,\cdots,w^*_N)$ can be
denoted by
$$W^*=e^{i\theta}(\mu,0,\cdots,0)U,\quad \mu\geqslant 0,$$
where $U$ is $(N, N)$ unitary matrix, therefore
$$
\overline{W^*}^tW^*=\overline{U}^t\left(\begin{array}{cc}\mu^2&0\\
0&\mathbf{0} \end{array}\right)U,$$
$$\quad M'_II^{(N)} +M''_I{\overline
{W^*}}^{t}W^{*}=\overline{U}^t\left(\begin{array}{cc}M'_I+M''_I\mu^2&0\\
0&M'_II^{(N-1)} \end{array}\right)U,
$$
but
$W^*\overline{W^*}^t=X=\mu^2.$ By calculations, one has
$$M'_I+M''_I\mu^2=M'_I+M''_IX=(N+1)Y^2+\frac{mn\lambda Y^2}
{(\lambda Y+K(m+n)+N-\lambda)^2} [(\lambda Y$$
$$+K(m+n)+N-\lambda)(Y-1)+Y(K(m+n)+N)]>0.$$

Therefore we proved that
$$\left(
\begin{array}{cc}
[\frac{1}{K}M'_IX+m+n+\frac{N}{K}]I&0\\
0&M'_I I+M_I''{\overline {W^*}}^t W^*
\end{array}\right)$$ is positive matrix. That is
$$\left(\frac{\partial^2\log[\det T_{\lambda I}(Z,W;\overline{Z},\overline{W})]}
{\partial z_i\partial
\overline{z}_j}\right)>0.$$
2.4. From 1.2.2., we know that
$$dzT_{\lambda I}(Z,W;\overline{Z},\overline{W})\overline{dz}^t=(\omega_{G_{\lambda}}(Y_I))^2$$
$$=dzJ_f|_{Z_0=Z}\left(
\begin{array}{cc}
[\frac{1}{K}\lambda YX+m+n+\frac{N}{K}]I&0\\
0& \lambda YI+\lambda Y^2{\overline {W^*}^t}W^*
\end{array}\right)\overline{J_f}^t|_{Z_0=Z}\overline{dz}^t
.$$ Where $Y, X$ are same as that in 1.1.2.

2.5. Let $dzJ_f|_{Z_0=Z}=(d\tilde{z}, d\Im)$, then similar from
1.2.4 to 1.2.5, one has
$$\begin{array}{lll}
(\omega_{G_{\lambda}}(Y_I))^2&=&(\frac{\lambda
YX}{K}+m+n+\frac{N}{K})|d\tilde{z}|^2+
(\lambda Y+\lambda Y^2\mu^2)|dv|^2+\lambda Y|d\tilde{W}|^2,\\
(\omega_{det}(Y_I))^2&=&(\frac{1}{K}M'_IX+m+n+\frac{N}{K})|d\tilde{z}|^2+(M'_I+M''_I\mu^2)|dv|^2
+M'_I|d\tilde{W}|^2.
\end{array}$$

2.6. Because $T_{\lambda I}(Z,W;\overline{Z},\overline{W})>0$ and
$(\frac{\partial^2\log[\det T_{\lambda
I}(Z,W;\overline{Z},\overline{W})]}{\partial z_i\partial
\overline{z}_j})>0,$ hence
$$\frac{1}{K}M'_IX+m+n+\frac{N}{K}>0,\quad M'_I+M''_I\mu^2 >0, \quad  M'_I >0,\quad Y >0.$$

2.7. Let
$$
\Phi_I(X)=\frac{\frac{1}{K}M'_IX+m+n+\frac{N}{K}}{\frac{\lambda
YX}{K}+m+n+\frac{N}{K}},\quad
\Psi_I(X)=\frac{M'_I+M''_I\mu^2}{\lambda Y+\lambda Y^2\mu^2},\quad
\Upsilon_I(X)=\frac{M'_I}{\lambda Y},
$$
then $\Phi_I(X),\Psi_I(X)$ and $\Upsilon_I(X)$ are the positive and
continuous functions of $X$ on $[0, 1).$ If the
$$\lim_{X\to 1}\Phi_I(X),\ \ \ \lim_{X\to 1}\Psi_I(X),
\ \ \ \lim_{X\to1}\Upsilon_I(X)$$ exist and positive, then
$\Phi_I(X), \Psi_I(X), \Upsilon_I(X)$ have the positive maximum
and positive minimum on $[0, 1)$.

2.8. We know the values of $M'_I, M''_I$ in 2.3., therefore one
can calculate the limits of $\Phi_I(X),\Psi_I(X)$ and
$\Upsilon_I(X)$ as $X\to 1$.

2.9. It is easy to show that
$$
\lim_{X\to1}\Phi_I(X)=\lim_{Y\to\infty}\Phi_I(X)=\frac{mn+N+1}{\lambda}.
$$
then there exists $0<\nu<\delta$ such that
$$0<\nu\leqslant\Phi(X)\leqslant \delta.$$

2.10. By the same method, one has
$$
\lim_{X\to
1}\Psi_I(X)=\lim_{Y\to\infty}\Psi_I(X)=\frac{mn+N+1}{\lambda}.
$$
and
$$
\lim_{X\to1}\Upsilon_I(X)=\lim_{Y\to\infty}\Upsilon_I(X)=\frac{mn+N+1}{\lambda}.
$$
Therefore there exist $\zeta,\ \eta,\ \rho$ and $\varrho$ such
that
$$
0<\zeta\leqslant\Psi(X)\leqslant\eta,\quad
0<\rho\leqslant\Upsilon(X)\leqslant\varrho.
$$

2.11. Let $a^2=\max\{\mu,\ \eta,\ \varrho\}$ and $b^2=\min\{\nu,\
\zeta,\ \rho\}$, then one has
$$
0<b\leqslant
\frac{\omega_{\det}(Y_I)}{\omega_{G_{\lambda}}(Y_I)}\leqslant a.
$$

2.12. Up to now, the following theorem is proved.

Theorem: The Ricci curvature of $\omega_{\lambda}(Y_I)$ on $Y_I$
is bounded from above and below by the negative constants, that is
$$-a\leqslant Ric_{\lambda I}=
-\frac{\omega_{det}(Y_I)}{\omega_{G_{\lambda}}(Y_I)}\leqslant
-b.$$

2.13. By using the same method, the Ricci curvature of
$\omega_{G_{\lambda II}}( \omega_{G_{\lambda III}},
\omega_{G_{\lambda IV}})$ on $Y_{II}(Y_{III}, Y_{IV})$ is also
bounded from above and below by the negative constants.

\section {Holomorphic sectional curvature of the new metrics}
In this section, the estimate of the holomorphic sectional
curvatures of new complete invariant metrics on Cartan-Hartogs
domains will be given. They are bounded form above and below by the
negative constants.

3.1. By the definition the holomorphic sectional curvature
$\omega_{\lambda I}(z, dz)$ of $\omega_{G_{\lambda}}(Y_I)$ on $Y_I$
has the following form:
$$\omega_{\lambda I}(z,
dz)=\frac{dz(-\overline{d}dT+dTT^{-1}\overline{dT}^t)
\overline{dz}^t}{(dzT\overline{dz}^t)^2},$$ where $T=T_{\lambda
I}(Z,W;\overline{Z},\overline{W}).$ Because the holomorphic
sectional curvature is invariant under the mapping of ${\rm
Aut}(Y_I)$. And for any $(Z,W)\in Y_I$, there exists a $f\in {\rm
Aut}(Y_I)$ such that $f(Z,W)=(0,W^*)$. Therefore it is sufficient
to compute the value of $\omega_{\lambda I}(z, dz)$ at point
$(0,W^*)$. In the following computation, $W$ stands for $W^*$ for
the sake of convenience. And if $Z^*=0$ then $|W^*|^2=X.$ Where
the $X$ can be found in 1.1.2.

3.2. Because
$$dT=\left(\begin{array}{cc}
dT_{11}& dT_{12}\\
dT_{21}& dT_{22}
\end{array}
\right),\ \ \ \overline{d}dT=\left(\begin{array}{cc}
\overline{d}dT_{11}& \overline{d}dT_{12}\\
\overline{d}dT_{21}& \overline{d}dT_{22}
\end{array}
\right).$$ Where the $T_{11}, T_{12}, T_{21}, T_{22}$ can be found
in 1.1.6.

3.3. We will compute the value  of $T,\ T^{-1},\ \overline{d}dT,\
dT,\ \overline{dT}$ at point $(0, W)$. Where
$$G(X)=Y^{\lambda},\ \log G(X)=M=\log Y^{\lambda},\
\frac{\partial \log G(X)}{\partial X}=M'=\lambda Y,$$
$$\frac{\partial^2
\log G(X)}{\partial X^2} =M''=\lambda Y^2.\ M'''=2\lambda Y^3, \
M^{(4)}=6\lambda Y^4.$$ Remark: Where the $M,\ M',\ M''$ are
different from that in 1.2.1 to 1.2.8.  By complicate
calculations, one has
$$T|_{Z=0}=\left(
\begin{array}{cc}
(\frac{1}{K}M'X+m+n+\frac{N}{K})I&0\\
0&M'I+M''{\overline W}^t W
\end{array}\right).$$
$$T^{-1}|_{Z=0}=
\left(\begin{array}{cc}
(\frac{1}{K}M'X+m+n+\frac{N}{K})^{-1}I&0\\
0&\frac{1}{M'}(I-(M'+M''X)^{-1}\overline{W}^tWM'')
\end{array}\right). $$
$$dT|_{Z=0}=\left(\begin{array}{cc}
\frac{XM''+M'}{K}\overline{W}dW^{t}I&0\\
\frac{XM''+M'}{K}\overline{W}^tdZ_1&M'''(\overline{W}dW^t)\overline{W}^tW+M''
(\overline{W}dW^{t}I+\overline{W}^tdW)
\end{array}\right).
$$
$$\begin{array}{ll}
\overline{d}dT_{11}|_{Z=0}&=\frac{X}{K^2}(XM''+M')|dZ_1|^{2}I
+(\frac{1}{K}M'X+m+n+\frac{N}{K})
(\overline{dZ}dZ^t\cdot\!\!\times I\\
\\
&+I \cdot\!\!\times \overline{dZ^t}dZ)+\frac{X}{K^2}(XM''+M')\overline{dZ_1^t}dZ_1\\
\\
&+\frac{1}{K}(XM'''+2M'')|W\overline{dW}^t|^{2}I+\frac{1}{K}(XM''+M')|dW|^2
I.
\end{array}$$
$$\begin{array}{ll}
\overline{d}dT_{12}|_{Z=0}=\frac{1}{K}(XM'''+2M'')(\overline{W}dW^t)(\overline{dZ_1^t}W)
+\frac{1}{K}(XM''+M')(\overline{dZ_1^t}dW). \end{array}$$
$$\begin{array}{ll}
\overline{d}dT_{21}|_{Z=0}=\frac{1}{K}(XM'''+2M'')(\overline{W}dW^t)(\overline
W^t{dZ_1}) +\frac{1}{K}(XM''+M')(\overline{dW}^tdZ_1).
\end{array}$$
$$
\begin{array}{ll}
\overline{d}dT_{22}|_{Z=0}&=\frac{1}{K}(XM''+M')|dZ_1|^2 I+
\frac{1}{K}(XM'''+2M'')|dZ_1|^2\overline{W}^tW\\
\\
&+(M''I+M'''\overline{W}^tW)|dW|^2
+M''\overline{dW}^tdW+M^{(4)}\overline{W}^tW
|W\overline{dW^t}|^2\\
\\
&+M'''[|W\overline{dW}^t|^2 I+({\overline W}dW^t) ({\overline
{dW}^t}W)+(W \overline{dW}^t)(\overline{W}^tdW)].
\end{array}$$

3.4. Let
$$-\overline{d}dT+dTT^{-1}\overline{dT'}|_{Z=0}
=\left(\begin{array}{cc}
R_{11}&R_{12}\\
R_{21}&R_{22} \end{array} \right),$$  then one has
$$
\begin{array}{lll}
R_{11}
&=&\left(\frac{1}{K^2}(XM''+M')^2 (\frac{1}{K}M'X+m+n+\frac{N}{K})^{-1}-\right.\\
\\
&&\left.\frac{1}{K}(XM'''+2M'')\right)I|\overline{W}dW^t|^2-
\frac{X}{K^2}(XM''+M')|dZ_1|^{2}I\\
\\
&&-\frac{1}{K}(XM''+M')|dW|^2I-\frac{X}{K^2}(XM''+M')\overline{dZ}^tdZ\\
\\
&&-(\frac{1}{K}M'X+m+n+\frac{N}{K})(\overline{dZ}dZ^t\cdot\!\!\times
I+
I\cdot\!\!\times \overline{dZ}^tdZ),\\
\\
R_{12}&=& [\frac{1}{K^2}(XM''+M')^2 (\frac{1}{K}M'X+m+n+\frac{N}{K})^{-1}\\
\\
&&-\frac{1}{K}(XM'''+2M'')]\overline{W}dW^t\overline{dZ_1}^tW-
\frac{1}{K}(XM''+M')(\overline{dZ_1}^tdW),\\
\\
R_{21}&=&\overline{R^t}_{12},\\
\\
R_{22}
&=&(M')^{-1}[M'''(XM'''+4M'')-(XM''+M')^{-1}M''(XM'''+2M'')^2]\\
\\
&&|W\overline{dW}^t|^2 \overline{W}^tW+[(M'')^2 (M')^{-1}-M''']
(\overline{W}dW^tI+\overline{W}^tdW)\\
\\
&&(\overline{dW}W^tI+\overline{dW}^tW)-\frac{1}{K}(XM''+M')|dZ_1|^2 I-M''|dW|^2 I\\
\\
&&+[\frac{1}{K^2}(XM''+M')^2 (\frac{1}{K}M'X+m+n+\frac{N}{K})^{-1}\\
\\
&&-\frac{1}{K}(XM'''+2M'')]|dZ_1|^2
\overline{W}^tW-M''\overline{dW}^t dW
-M^{(4)}|W\overline{dW}^t|^2 \overline{W}^tW.\\
\\
\end{array}
$$

3.5. Because
$$
\begin{array}{ll}
\ \ \ \ dz\left(-\overline{d}dT+dTT^{-1}\overline{dT}^t\right)\overline{dz}^t|_{Z=0}\\
\\
=(dZ_1,dW)\left(\begin{array}{cc}
R_{11}&R_{12}\\
R_{21}&R_{22}\end{array}
\right)(\overline{dZ}_1,\overline{dW}^t)\\
\\
=dZ_1R_{11}\overline{dZ}_1^t+dWR_{21}\overline{dZ}_1^t+dZ_1R_{12}\overline{dW}^t
+dWR_{22}\overline{dW}^t.\end{array}$$
By calculations one has
$$\begin{array}{ll}
dz[-\overline{d}dT+dTT^{-1}\overline{dT}^t]\overline{dz}^t|_{Z=0}
\end{array}$$
$$\begin{array}{ll}
=P_1|W\overline{dW}^t|^4+P_{12}|W\overline{dW}^t|^2 |dW|^2+P_2
|dW|^4\\
\\
\ \ \
+Q_1 |dW|^2 |dZ_1|^2 +Q_2|W\overline{dW}^t|^2 |dZ_1|^2 \\
\\
\ \ \
+R|dZ_1|^4-(\frac{1}{K}M'X+m+n+\frac{N}{K})dZ_1(\overline{dZ}dZ^t\cdot
\!\!\times I+I\cdot\!\!\times \overline{dZ^t}dZ)\overline{dZ}_1^t.\\
\\
\end{array}
$$
Where
$$
\begin{array}{ll}
P_1=\frac{1}{M'}[M'''(XM'''+4M'')-(XM''+M')^{-1}M''(XM'''+2M'')^2]-M^{(4)}\\
\\
\ \quad=-2\lambda Y^4,\\
\\
P_{12}=4[(M'')^2(M')^{-1}-M''']=-4\lambda Y^3,\\
\\
P_2=-2M''=-2\lambda Y^2,\\
\\
Q_1=-\frac{4}{K}(XM''+M')=-\frac{4\lambda Y^2}{K},\\
\\
Q_2=4[\frac{1}{K^2}(XM''+M')^2(\frac{1}{K}M'X+m+n+\frac{N}{K})^{-1}-
\frac{1}{K}(XM'''+2M'')]\\
\\
\ \quad=\frac{4}{K^2}\lambda^2Y^4(\frac{\lambda
Y}{K}+m+n+\frac{N-\lambda}{K})^{-1}-\frac{8}{K}\lambda Y^3\\
\\
\ \quad=\frac{4}{K}\lambda^2Y^4(\lambda Y+M_1)^{-1} -\frac{8}{K}\lambda Y^3,\\
\\
R=-\frac{2}{K^2}(XM''+M')X=-\frac{2}{K^2}\lambda (Y^2-Y).\\
\\
\end{array}$$ where $M_1=(m+n)K+N-\lambda.$

3.6. By calculations one
has$$dZ_1(\overline{dZ}dZ^t\cdot\!\!\times I+ I\cdot\!\!\times
\overline{dZ}^tdZ)\overline{dZ_1}^t =2{\rm
tr}(dZ\overline{dZ}^tdZ\overline{dZ}^t).$$ Let $$\Omega_1=
P_1|W\overline{dW}^t|^4+P_{12}|W\overline{dW}^t|^2 |dW|^2 +P_2
|dW|^4+Q_1 |dW|^2 |dZ_1|^2 $$ $$+Q_2|W\overline{dW}^t|^2
|dZ_1|^2+R|dZ_1|^4- 2K^{-1}(\lambda Y+M_1){\rm
tr}(dZ\overline{dZ}^tdZ\overline{dZ}^t),$$ and
$$\Omega_2=[K^{-1}(\lambda Y+M_1)|dZ_1|^2+\lambda Y|dW|^2+\lambda
Y^2|W\overline{dW}^t|^2)]^2.$$ Hence
$$
\omega_{\lambda I}(z, dz)|_{Z=0}=\frac{\Omega_1}{\Omega_2}.$$

3.7. Let $Z$ be the $(m,n)$ complex matrix, then it is easy to get
$${\rm tr}(Z\overline{Z}^tZ\overline{Z}^t)\leqslant {\rm tr}(Z\overline{Z}^t)
{\rm tr}(Z\overline{Z}^t) \leqslant m{\rm
tr}(Z\overline{Z}^tZ\overline{Z}^t).
$$ Therefore one has
$$\frac{2}{mK}(\lambda Y+M_1)|dZ_1|^4\leqslant \frac{2}{K}(\lambda Y+M_1)
{\rm tr}(dZ\overline{dZ}^tdZ\overline{dZ}^t) \leqslant
\frac{2}{K}(\lambda Y+M_1)|dZ_1|^4.$$

3.8. The lower bound of holomorphic sectional curvature

3.8.1. From 3.6., we know that the $P_1$ and $\lambda^2 Y^4$ are
the coefficients of $|W\overline{dW}^t|^4$ in $\Omega_1$ and
$\Omega_2$ respectively, let
$$\Phi_1=\frac{P_1}{\lambda^2 Y^4},$$ the
$P_{12}$ and $2\lambda^2 Y^3$ are the coefficients of
$|W\overline{dW}^t|^2 |dW|^2$ in $\Omega_1$ and $\Omega_2$
respectively, let
$$\Phi_2=\frac{P_{12}}{2\lambda^2 Y^3},$$ the $P_2$ and
$\lambda^2 Y^2$ are the coefficients of $|dW|^4$ in $\Omega_1$ and
$\Omega_2$ respectively, let
$$\Phi_3=\frac{P_2}{\lambda^2 Y^2},
$$ the $Q_1$ and $2K^{-1}\lambda Y(\lambda Y+M_1)$
are the coefficients of $|dW|^2 |dZ_1|^2$ in $\Omega_1$ and
$\Omega_2$ respectively, let
$$\Phi_4=\frac{Q_1}{2K^{-1}\lambda Y(\lambda Y+M_1)},$$
the $Q_2$ and $2K^{-1}(\lambda Y+M_1)\lambda Y^2$ are the
coefficients of $|W\overline{dW}^t|^2 |dZ_1|^2$ in $\Omega_1$ and
$\Omega_2$ respectively, let
$$\Phi_5=\frac{Q_2}{2K^{-1}(\lambda Y+M_1)\lambda Y^2}$$ and let
$$
\Phi_6=\frac{R|dZ_1|^4-2K^{-1}(\lambda Y+M_1) {\rm
tr}(dZ\overline{dZ}^tdZ\overline{dZ}^t)}{K^{-2}(\lambda
Y+M_1)^2|dZ_1|^4}.
$$
Then one has
$$\Phi_1=-\frac{2}{\lambda},\ \ \Phi_2=-\frac{2}{\lambda},\ \
\Phi_3=-\frac{2}{\lambda},\ \ \Phi_4=-\frac{2Y}{\lambda Y+M_1}.$$

3.8.2. Because
$$
\Phi_5=-\frac{2Y(\lambda Y+2M_1)}{(\lambda Y+M_1)^2}\geqslant
-\frac{2Y[\lambda Y+2(m+n)K+2N]}{(\lambda Y+M_1)^2}:=\Phi_{51}.
$$
And by using the inequality in 3.7., one has
$$
\Phi_6\geqslant
\frac{R|dZ_1|^4-2K^{-1}(\lambda Y+M_1)|dZ_1|^4}{K^{-2}(\lambda
Y+M_1)^2|dZ_1|^4}=\frac{-2\lambda (Y^2-Y)-2K(\lambda
Y+M_1)}{(\lambda Y+M_1)^2}:=\Phi_{61}.
$$

3.8.3. It is easy to see that the $\Phi_1,\ \Phi_2,\ \Phi_3,\
\Phi_4,\ \Phi_{51},\ \Phi_{61}$ are the negative continues
functions of $Y$ on the interval $[1, \infty)$. If $Y\to\infty$,
then their limits are existent and are the negative numbers. Hence
all of $\Phi_1,\ \Phi_2,\ \Phi_3,\ \Phi_4,\ \Phi_{51},\ \Phi_{61}$
have the negative minimums on $[1, \infty)$ respectively. Let $-a$
be the smallest one of them. Then it is easy to show that
$\Omega_1\geqslant -a\Omega_2$, that is
$$\omega_{\lambda I}(z, dz)|_{Z=0}\geqslant -a.$$

3.9. The upper bound of holomorphic sectional curvature

The $\omega_{\lambda I}(z, dz)|_{Z=0}$ can be rewritten as
$$\omega_{\lambda I}(z, dz)|_{Z=0}=-C+\frac{\Omega_3}{\Omega_4}, \ \ C>0.$$
Where
$$
\Omega_3=P^*_1|W\overline{dW}^t|^4+P^*_{12}|W\overline{dW}^t|^2
|dW|^2 +P^*_2 |dW|^4+Q^*_1 |dW|^2 |dZ_1|^2
$$
$$
+Q^*_2|W\overline{dW}^t|^2 |dZ_1|^2+R^*|dZ_1|^4- 2K^{-1}(\lambda
Y+M_1){\rm tr}(dZ\overline{dZ}^tdZ\overline{dZ}^t),
$$
$$
\Omega_4=[K^{-1}(\lambda Y+M_1)|dZ_1|^2+\lambda Y|dW|^2+\lambda
Y^2|W\overline{dW}^t|^2)]^2,
$$
$$P_1^*=P_1+C\lambda^2Y^4=-\lambda Y^4(2-C\lambda),$$
$$P_{12}^*=P_{12}+aC\lambda^2Y^3=-2\lambda Y^3(2-C\lambda),\
P_2^*=P_2+C\lambda^2Y^2=-\lambda Y^2(2-C\lambda),$$
$$Q_1^*=Q_1+2CK^{-1}\lambda Y(\lambda Y+M_1)=-4K^{-1}\lambda
Y^2+2CK^{-1}\lambda Y(\lambda Y+M_1),$$
$$Q^*_2=Q_2+2CK^{-1}\lambda
Y^2(\lambda Y+M_1)$$ $$=4K^{-1}\lambda^2Y^4(\lambda Y+M_1)^{-1}
-8K^{-1}\lambda Y^3+2CK^{-1}\lambda Y^2(\lambda Y+M_1),$$
$$R^*=R+CK^{-2}(\lambda Y+M_1)^2=-2\lambda
K^{-2}(Y^2-Y)+CK^{-2}(\lambda Y+M_1)^2.$$ If
$$P^*_1\leqslant 0,\ \ \ P_{12}^*\leqslant 0,\ \ \ P^*_2\leqslant 0,$$
and
$$Q^*_1 |dW|^2 |dZ_1|^2
+Q^*_2|W\overline{dW}^t|^2 |dZ_1|^2\leqslant 0,$$
$$R^*|dZ_1|^4-2K^{-1}(\lambda Y+M_1) {\rm
tr}(dZ\overline{dZ}^tdZ\overline{dZ}^t)\leqslant 0,$$ then
$$\omega_{\lambda I}(z, dz)|_{Z=0}\leqslant -C.$$

3.9.1. It is easy to see that if $C\leqslant \frac{2}{\lambda}$,
then $P_1^*\leqslant 0,\ P_{12}^*\leqslant 0,\ P_2^*\leqslant 0.$

3.9.2. Because
$|W\overline{dW^t}|^2\leqslant|W|^2|dW|^2=X|dW|^2=(1-Y^{-1})|dW|^2$,
and if $C\leqslant \frac{2\lambda Y}{\lambda Y+M_1}:=\Phi_{42}, $
then $Q^*_1\leqslant 0$. At this time one has $$
Q_1^*|dW|^2|dZ_1|^2+Q_2^*|W\overline{dW^t}|^2|dZ_1|^2\leqslant
(Q_1^*X^{-1}+Q_2^*)|W\overline{dW^t}|^2|dZ_1|^2.$$ Then by
calculations, one has
$$
Q_1^*X^{-1}+Q_2^*=\frac{2C\lambda
Y^3(\lambda Y+M_1)}{K(Y-1)}-\frac{4\lambda
Y^3[\lambda(Y-1)^2+(M_1+\lambda)(2Y-1)]}{K(Y-1)(\lambda Y+M_1)}.
$$
Therefore if
$$C\leqslant \frac{2[\lambda(Y-1)^2+(M_1+\lambda)(2Y-1)]}{(\lambda Y+M_1)^2}:=\Phi_{52},$$
and $$C\leqslant\frac{2Y}{(\lambda Y+M_1)}:=\Phi_{42},$$ where
$M_1=(m+n)K+N-\lambda,$ then
$$Q_1^*|dW|^2|dZ_1|^2+Q_2^*|W\overline{dW^t}|^2|dZ_1|^2\leqslant 0.$$

3.9.3. By using the inequality in 3.7., one has
$$
R^*|dZ_1|^4- \frac{2(\lambda Y+M_1){\rm
tr}(dZ\overline{dZ}^tdZ\overline{dZ}^t)}{K}\leqslant R^*|dZ_1|^4-
\frac{2(\lambda Y+M_1)}{mK}|dZ_1|^4.
$$
By calculations, if
$$
C\leqslant \frac{2[\lambda Y^2+Km^{-1}(\lambda Y+M_1)]}{(\lambda
Y+M_1)^2}:=\Phi_{62},
$$
then
$$
R^*|dZ_1|^4- 2K^{-1}(\lambda Y+M_1) {\rm
tr}(dZ\overline{dZ}^tdZ\overline{dZ}^t)\leqslant 0.
$$

3.9.4. Because $\Phi_{42},\ \Phi_{52},\ \Phi_{62}$ are the
positive continues functions of $Y$ on the interval $[1, \infty)$.
It is easy to show that, when $Y\to \infty$, the limits of $
\Phi_{42},\ \Phi_{52},\ \Phi_{62}$ are existent and are equal to
the positive numbers $\frac{2}{\lambda}$. Then $\Phi_{42},\
\Phi_{52},\ \Phi_{62}$ have the positive minimums  on $[1,
\infty)$ respectively. Let $b$ be the smallest one of them. Then
if $C\leqslant b$, one has
$$Q_1^*|dW|^2|dZ_1|^2+Q_2^*|W\overline{dW^t}|^2|dZ_1|^2\leqslant 0,$$
$$R^*|dZ_1|^4- 2K^{-1}(\lambda Y+M_1)
{\rm tr}(dZ\overline{dZ}^tdZ\overline{dZ}^t)\leqslant 0.$$

3.9.5. From 3.9.1 to 3.9.4, there exists
$$C\leqslant\min\{b,\frac{2}{\lambda}\},$$
and $C>0$ such that
$$\omega_{\lambda I}(z, dz)\leqslant -C.$$

By the 3.8.3 and 3.9.5., one has the following theorem.

3.10. Theorem: There exists positive constant $a$, $C$ dependent
on $Y_I$, $\lambda$ such that the holomorphic sectional curvature
$\omega_{\lambda I}(z, dz)$ of metric $\omega_{G_{\lambda}}(Y_I)$
on $Y_I$ satisfies
$$-a\leqslant\omega_{\lambda I}(z, dz)\leqslant -C.$$

This theorem is also true for the other Cartan-Hartogs domains.

\section {Bergman metric is equivalent to the Einstein-K\"ahler metric}
We proved that the new complete invariant metrics are equivalent
to the Bergman metric on Cartan-Hartogs domains. We will prove
that these new metrics are also equivalent to the
Einstein-K\"ahler metric on Cartan-Hartogs. Therefore the Bergman
metric is equivalent to the Einstein-K\"ahler metric on
Cartan-Hartogs domain.  By using the Yau's Schwarz lemma and the
theorem 3.10, we can prove that the Bergman metric is equivalent
to the Einstein-K\"ahler metric on Cartan-Hartogs domain.

4.1. Yau's Schwarz lemma[33]: Let $f : (M^m, g)\rightarrow (N^n,
h)$ be a holomorphic map between K\"ahler manifolds where $M$ is
complete and $Ric(g)\geqslant -cg$ with $c\geqslant 0$.

(1) if the holomorphic sectional curvature of $N$ is bounded above
by a negative constant, then $f^*h\leqslant \tilde{c}g$ for some
constant $\tilde{c}$.

(2) If $m=n$ and the Ricci curvature of $N$ is bounded above by a
negative constant, then $f^*\omega_h^n\leqslant
\tilde{c}\omega^n_g$ for some constant $\tilde{c}$.

Where $\omega^n_h$ and $\omega^n_g$ are the volume element for
$(M^m, g)$ and $(N^n, h)$ respectively.

4.2. Consider the identity map
$$id: (Y_I, \omega_{EK}(Y_I))\rightarrow (Y_I, \omega_{G_{\lambda}}(Y_I)), $$
Because the holomorphic sectional curvature of
$\omega_{G_{\lambda}}(Y_I)$ is bounded above by a negative constant
. Yau's Schwarz lemma (1) implies
$$\omega_{G_{\lambda}}(Y_I)\leqslant C_1\omega_{EK}(Y_I).$$
Consider the identity map again
$$id: (Y_I, \omega_{G_{\lambda}}(Y_I))\rightarrow (Y_I, \omega_{EK}(Y_I)).
$$
Because the Ricci curvature of $\omega_{G_{\lambda}}(Y_I)$ is
bounded below by a negative constant. Yau's Schwarz lemma (2)
implies
$$\omega^{mn+N}_{EK}(Y_I)\leqslant
C_0\omega^{mn+N}_{G_{\lambda}}(Y_I).$$Where $mn+N$ is the
dimension of $Y_I.$ This inequality implies
$$\det[T_{EK I}(Z,W;\overline{Z},\overline{W})]
\leqslant C_0\det[T_{\lambda I}(Z,W;\overline{Z},\overline{W})].$$
Because $T_{\lambda I}(Z,W;\overline{Z},\overline{W})>0,\ T_{EK
I}(Z,W;\overline{Z},\overline{W})>0$. Then from the following
proposition, one has $$\omega_{EK}(Y_I)\leq
C_2\omega_{G_{\lambda}}(Y_I).$$

Proposition: Let $A$ and $B$ be positive definite $n\times n$
Hermitian matrices and let $\alpha,\ \beta$ be positive constants
such that $B\geqslant \alpha A$ and $\det(B)\leqslant
\beta\det(A)$. Then there is a constant $\gamma >0$ depending on
$\alpha,\ \beta$ and $n$ such that $B\leqslant \gamma A$.

Up to now we proved that

Theorem: The Bergman metric is equivalent to the
Einstein-K\"ahler metric on $Y_I$.

This theorem are also true for the other Cartan-Hartogs domains.
Thus the Yau's conjecture is true for the Cartan-Hartogs domains.

4.3. If $Y_I$ is convex, then $\omega_B(Y_I),\ \omega_C(Y_I),\ \
\omega_K(Y_I),\ \omega_{EK}(Y_I),\ \omega_{G_{\lambda}}(Y_I)$ are
equivalence on $Y_I$. This fact is also true for the other convex
Cartan-Hartogs domains.

4.4. Because  the holomorphic sectional curvature of
$\omega_{G_{\lambda}}(Y_I)$ is bounded above by negative constant,
then by the ref.[41, p.136], one has
$$\omega_{G_{\lambda}}(Y_I)\leqslant \beta \omega_K(Y_I) $$ Hence
$$\omega_B(Y_I)\leqslant \beta_1
\omega_K(Y_I)$$ and
$$\omega_{EK}(Y_I)\leqslant \beta_2 \omega_K(Y_I).$$

4.5. Because $ \omega_C(Y_I)\leqslant 2\omega_B(Y_I),$ then
$$\omega_C(Y_I)\leqslant \beta_3 \omega_{EK}(Y_I).$$

Where $\beta,\ \beta_1,\ \beta_2,\ \beta_3 $ are the positive
constants. The facts in 4.4 and 4.5 are also true for the other
Cartan-Hartogs domains.

\end{document}